\documentclass[a4paper,10pt]{article}

\usepackage{latexsym}
\usepackage{amsmath,amsthm,amsfonts}
\usepackage{graphicx}

\newcommand{\bN}{\mathbb{N}}

\newcommand{\bZ}{\mathbb{Z}}

\theoremstyle{plain}

\newtheorem{Theorem}{Theorem}

\newtheorem{Proposition}[Theorem]{Proposition}

\theoremstyle{remark}

\begin{document}

\title{On the Number of Hamiltonian Groups}

\author{
  Boris Horvat\thanks{\texttt{Boris.Horvat@fmf.uni-lj.si}} \\
  University of Ljubljana, Slovenia\\
  Ga\v sper Jakli\v c\thanks{\texttt{Gasper.Jaklic@fmf.uni-lj.si}} \\
  University of Ljubljana, Slovenia \\
  Toma\v{z} Pisanski\thanks{\texttt{Tomaz.Pisanski@fmf.uni-lj.si}} \\
  University of Ljubljana and University of Primorska, Slovenia
}
\date{\today}

\maketitle

\begin{abstract}
Finite hamiltonian groups are counted. The sequence of numbers of all groups of order $n$ all
whose subgroups are normal and the sequence of numbers of all groups of order less or equal to 
$n$ all whose subgroups are normal are presented.
\end{abstract}


\section{Introduction}
Subgroups of abelian groups are abelian and hence self-conjugate or \emph{normal}. A nonabelian 
group all of whose subgroups are normal is called \emph{hamiltonian} \cite{C56,We04h}. 
Let $\cal A$ denote the class of abelian groups and let $\cal H$ denote the class of hamiltonian
groups. In topological graph theory \cite{GT87,Wh01}, hamiltonian groups have been studied in
the past \cite{PT89,PW88,PT89-1}. For several classes of hamiltonian groups the genus is known
exactly. For abelian and hamiltonian groups, there are structural theorems available. 
We note in passing that here we use a different structure theorem. For instance, the cyclic 
group $\bZ_{15}$ can be written as $\bZ_3 \times \bZ_5.$ Since it can be generated by a single 
generator, the former form is preferred in the topological graph theory over the latter.  
In this paper we determine the number $h(n)$ of hamiltonian groups of order $n$ and the number
$b(n)$ of all groups of order $n$ with the property, that all their subgroups are normal. We also 
determine the number $v(n)$ of all hamiltonian groups of order $\le n$ and the number $w(n)$ of 
all groups of order $\le n$ with the property, that all their subgroups are normal. 

\section{Results}
Before we study hamiltonian groups we will recall the structure of finite abelian groups 
\cite{We04a}. 
Let $\pi(m)$ denote a partition of a natural number $m$, where
$$\pi(m) := \{m_1, m_2, \dots, m_s\},  
$$
such that $m = \sum_{k=1}^s m_k$ and $m_i \geq m_j$ for all $1 \leq i < j \leq s$.
For $c \in \bN$ let $c^{\pi(m)} := \{c^{m_1}, c^{m_2}, \dots, c^{m_s}\}$ and let
$A(n_1, n_2, \dots, n_r)$ denote the direct product of cyclic groups
$$A(n_1,n_2, \dots, n_r) := \bZ_{n_1} \times \bZ_{n_2} \times \dots \times
\bZ_{n_r}.
$$
Let $G$ be a finite abelian group of order $n$. Let us write down the prime decomposition of
$n$ as
$$n = \prod_{k=1}^\ell p_k^{\alpha_k}.$$
It is well-known that $G$ is isomorphic to the group
$$
G \approx A\left(p_1^{\pi(\alpha_1)},p_2^{\pi(\alpha_2)}, \dots, p_\ell^{\pi(\alpha_\ell)}\right).
$$
Let $a(n)$ denote the number of abelian groups of order $n$ and let $P(n)$ denote the number of
partitions of the integer $n.$ The previous discussion gives a proof to the following result.

\begin{Proposition}
The number $a(n)$ of abelian groups of order $n$ is given by
$\prod_{i=1}^\ell P(\alpha_i)$ where $n = \prod_{k=1}^\ell p_k^{\alpha_k}$ is the prime
decomposition of $n$.
\end{Proposition}

\noindent The initial 200 values of the sequence $a(n)$ are given in Table 
\ref{table-initvalues-an}.

\begin{table}[h]
\begin{tiny}
\begin{center}
\begin{tabular}{|c|ccccc ccccc ccccc ccccc|}
\hline \hline
$n$&1& & & & 5& & & & & 10& & & & & 15& & & & & 20\\
\hline
0&   1 & 1 & 1 & 2 & 1 & 1 & 1 & 3 & 2 & 1 & 1 & 2 & 1 & 1 & 1 & 5 & 1 & 2 & 1 & 2 \\
20 &   1 & 1 & 1 & 3 & 2 & 1 & 3 & 2 & 1 & 1 & 1 & 7 & 1 & 1 & 1 & 4 & 1 & 1 & 1 & 3 \\ 
40 &   1 & 1 & 1 & 2 & 2 & 1 & 1 & 5 & 2 & 2 & 1 & 2 & 1 & 3 & 1 & 3 & 1 & 1 & 1 & 2 \\
60 &   1 & 1 & 2 & 11 & 1 & 1 & 1 & 2 & 1 & 1 & 1 & 6 & 1 & 1 & 2 & 2 & 1 & 1 & 1 & 5 \\ 
80 &   5 & 1 & 1 & 2 & 1 & 1 & 1 & 3 & 1 & 2 & 1 & 2 & 1 & 1 & 1 & 7 & 1 & 2 & 2 & 4 \\
100 & 1 & 1 & 1 & 3 & 1 & 1 & 1 & 6 & 1 & 1 & 1 & 5 & 1 & 1 & 1 & 2 & 2 & 1 & 1 & 3 \\
120 & 2 & 1 & 1 & 2 & 3 & 2 & 1 & 15 & 1 & 1 & 1 & 2 & 1 & 1 & 3 & 3 & 1 & 1 & 1 & 2 \\
140 & 1 & 1 & 1 & 10 & 1 & 1 & 2 & 2 & 1 & 2 & 1 & 3 & 2 & 1 & 1 & 2 & 1 & 1 & 1 & 7 \\
160 & 1 & 5 & 1 & 2 & 1 & 1 & 1 & 3 & 2 & 1 & 2 & 2 & 1 & 1 & 2 & 5 & 1 & 1 & 1 & 4 \\
180 & 1 & 1 & 1 & 3 & 1 & 1 & 1 & 2 & 3 & 1 & 1 & 11 & 1 & 1 & 1 & 4 & 1 & 2 & 1 & 6 \\
\hline \hline
\end{tabular}
\caption{The initial values of $a(n), \, n=1, 2, \dots, 200$, (\cite{SloaneA000688}, A000688).}
\label{table-initvalues-an}
\end{center}
\end{tiny}
\end{table}

\noindent Note that the sequence $\left\{a(n)\right\}_{n \in \bN}$ can not contain multiples
of primes in the sequence $s:=\{13, 17, 19, 23, 29, 31, 37, \ldots\}$ since 
$P(n) \ne k \cdot s_i, \forall n, i, k \in \bN$, (see \cite{SloaneA046064}). 
The number $a(n)$ depends only on the prime signature of $n$. For example, both $24 = 2^3 \cdot 3^1$ and $375 = 3^1 \cdot 5^3$ have the prime signature $(3, 1)$, therefore $a(375) = a(24)$.

A similar structural theorem holds for hamiltonian groups. A hamiltonian group $H$ is
isomorphic to a direct product of the quaternion group $Q$ of order $8$, an elementary
abelian group $E$ of exponent $2$ and an abelian group $A$ of odd order
$$
	H \approx Q \times E \times A \approx Q \times \bZ_{2^k} \times A,
$$
where $|Q|=8=2^3$, $|E|=2^k$ and $|A| \ne 0\; (\mbox{mod}\; 2)$. Therefore $|H|=2^{3+k} |A|$.
Let $n$ be an arbitrary natural number. We can write $n$ uniquely in the form 
$n = 2^e \cdot o$ where $e = e(n) \geq 0$ and $o = o(n)$ is an odd number. These
results give the number of hamiltonian groups of order $n$.

\begin{Proposition}
Let $n = 2^e \cdot o$, where $e = e(n) \geq 0$ and $o = o(n)$ is an odd number. The number 
$h(n)$ of hamiltonian groups of order $n$ is given by
\begin{displaymath}
h(n) = \left\{ \begin{array} {ll} 
0, & e(n) < 3; \\
a\left(o(n)\right), &  \textrm{otherwise}.
\end{array} \right.
\end{displaymath}
\end{Proposition}

\noindent The initial 200 values of the sequence $h(n)$ are given in Table 
\ref{table-initvalues-hn}.

\begin{table}[h]
\begin{tiny}
\begin{center}
\begin{tabular}{|c|ccccc ccccc ccccc ccccc|}
\hline \hline
$n$&1& & & & 5& & & & & 10& & & & & 15& & & & & 20\\
\hline
0 & 0 & 0 & 0 & 0 & 0 & 0 & 0 & 1 & 0 & 0 & 0 & 0 & 0 & 0 & 0 & 1 & 0 & 0 & 0 & 0 \\
20 & 0 & 0 & 0 & 1 & 0 & 0 & 0 & 0 & 0 & 0 & 0 & 1 & 0 & 0 & 0 & 0 & 0 & 0 & 0 & 1 \\
40 & 0 & 0 & 0 & 0 & 0 & 0 & 0 & 1 & 0 & 0 & 0 & 0 & 0 & 0 & 0 & 1 & 0 & 0 & 0 & 0 \\
60 & 0 & 0 & 0 & 1 & 0 & 0 & 0 & 0 & 0 & 0 & 0 & 2 & 0 & 0 & 0 & 0 & 0 & 0 & 0 & 1 \\
80 & 0 & 0 & 0 & 0 & 0 & 0 & 0 & 1 & 0 & 0 & 0 & 0 & 0 & 0 & 0 & 1 & 0 & 0 & 0 & 0 \\
100 & 0 & 0 & 0 & 1 & 0 & 0 & 0 & 0 & 0 & 0 & 0 & 1 & 0 & 0 & 0 & 0 & 0 & 0 & 0 & 1 \\
120 & 0 & 0 & 0 & 0 & 0 & 0 & 0 & 1 & 0 & 0 & 0 & 0 & 0 & 0 & 0 & 1 & 0 & 0 & 0 & 0 \\
140 & 0 & 0 & 0 & 2 & 0 & 0 & 0 & 0 & 0 & 0 & 0 & 1 & 0 & 0 & 0 & 0 & 0 & 0 & 0 & 1 \\
160 & 0 & 0 & 0 & 0 & 0 & 0 & 0 & 1 & 0 & 0 & 0 & 0 & 0 & 0 & 0 & 1 & 0 & 0 & 0 & 0 \\
180 & 0 & 0 & 0 & 1 & 0 & 0 & 0 & 0 & 0 & 0 & 0 & 1 & 0 & 0 & 0 & 0 & 0 & 0 & 0 & 2 \\
\hline \hline
\end{tabular}
\caption{The initial values of $h(n), \, n=1, 2, \dots, 200$.}
\label{table-initvalues-hn}
\end{center}
\end{tiny}
\end{table}

Combining abelian and hamiltonian groups of order $n$ we may now
give the number $b(n) := a(n) + h(n)$ of all groups of order $n$ all of whose
subgroups are normal. The initial 300 values of the sequence $b(n)$ are given
in Table \ref{table-initvalues-bn}.

\begin{table}[h]
\begin{tiny}
\begin{center}
\begin{tabular}{|c|ccccc ccccc ccccc ccccc|}
\hline \hline
$n$&1& & & & 5& & & & & 10& & & & & 15& & & & & 20\\
\hline
0&    1& 1& 1& 2& 1&   1& 1& 4& 2& 1&   1& 2& 1& 1& 1&   6& 1& 2& 1& 2\\
20&   1& 1& 1& 4& 2&   1& 3& 2& 1& 1&   1& 8& 1& 1& 1&   4& 1& 1& 1& 4\\
40&   1& 1& 1& 2& 2&   1& 1& 6& 2& 2&   1& 2& 1& 3& 1&   4& 1& 1& 1& 2\\
60&   1& 1& 2& 12& 1&  1& 1& 2& 1& 1&   1& 8& 1& 1& 2&   2& 1& 1& 1& 6\\
80&   5& 1& 1& 2& 1&   1& 1& 4& 1& 2&   1& 2& 1& 1& 1&   8& 1& 2& 2& 4\\
100&   1& 1& 1& 4& 1&   1& 1& 6& 1& 1&   1& 6& 1& 1& 1&   2& 2& 1& 1& 4\\
120&   2& 1& 1& 2& 3&   2& 1& 16& 1& 1&  1& 2& 1& 1& 3&   4& 1& 1& 1& 2\\
140&   1& 1& 1& 12& 1&  1& 2& 2& 1& 2&   1& 4& 2& 1& 1&   2& 1& 1& 1& 8\\
160&   1& 5& 1& 2& 1&   1& 1& 4& 2& 1&   2& 2& 1& 1& 2&   6& 1& 1& 1& 4\\
180&   1& 1& 1& 4& 1&   1& 1& 2& 3& 1&   1& 12& 1& 1& 1&  4& 1& 2& 1& 8\\
200&   1& 1& 1& 2& 1&   1& 2& 6& 1& 1&   1& 2& 1& 1& 1&  12& 1& 1& 1& 2\\
220&   1& 1& 1& 8& 4&   1& 1& 2& 1& 1&   1& 4& 1& 2& 1&   2& 1& 1& 1& 6\\
240&   1& 2& 7& 2& 2&   1& 1& 4& 1& 3&   1& 4& 1& 1& 1&  23& 1& 1& 1& 2\\
260&   2& 1& 1& 4& 1&   1& 1& 2& 1& 3&   1& 6& 1& 1& 2&   2& 1& 1& 2& 4\\
280&   1& 1& 1& 2& 1&   1& 1& 16& 2& 1&  1& 2& 1& 2& 1&   4& 3& 1& 1& 4\\
\hline \hline
\end{tabular}
\caption{The initial values of $b(n), \, n=1,2,\dots,300$.}
\label{table-initvalues-bn}
\end{center}
\end{tiny}
\end{table}

The number $u(n)$ of all abelian groups of order $\le n$ is presented in 
\cite{SloaneA063966}. The initial 100 values of the sequence $u(n)$ are
given in Table \ref{table-initvalues-un}.

\begin{table}[h]
\begin{tiny}
\begin{center}
\begin{tabular}{|c|ccccc ccccc|}
\hline \hline
$n$&1& & & & 5& & & & & 10\\
\hline 
0 & 1 & 2 & 3 & 5 & 6 & 7 & 8 & 11 & 13 & 14\\
10 & 15 & 17 & 18 & 19 & 20 & 25 & 26 & 28 & 29 & 31\\
20 & 32 & 33 & 34 & 37 & 39 & 40 & 43 & 45 & 46 & 47\\
30 & 48 & 55 & 56 & 57 & 58 & 62 & 63 & 64 & 65 & 68\\
40 & 69 & 70 & 71 & 73 & 75 & 76 & 77 & 82 & 84 & 86\\
50 & 87 & 89 & 90 & 93 & 94 & 97 & 98 & 99 & 100 & 102\\
60 & 103 & 104 & 106 & 117 & 118 & 119 & 120 & 122 & 123 & 124\\
70 & 125 & 131 & 132 & 133 & 135 & 137 & 138 & 139 & 140 & 145\\
80 & 150 & 151 & 152 & 154 & 155 & 156 & 157 & 160 & 161 & 163\\
90 & 164 & 166 & 167 & 168 & 169 & 176 & 177 & 179 & 181 & 185\\
\hline \hline
\end{tabular}
\caption{The initial values of $u(n), \, n=1, 2, \dots, 100$, (\cite{SloaneA063966}, A063966).}
\label{table-initvalues-un}
\end{center}
\end{tiny}
\end{table}

Let $v(n)$ be the number of all hamiltonian groups of order $\le n$ and let $w(n)$ be the
number of all groups of order $\le n$ all of whose subgroups are normal. The initial 200
values of the sequences $v(n)$ and $w(n)$ are given in 
Table \ref{table-initvalues-vn} and Table \ref{table-initvalues-wn}, respectively.

\begin{table}[h]
\begin{tiny}
\begin{center}
\begin{tabular}{|c|ccccc ccccc|}
\hline \hline
$n$&1& & & & 5& & & & & 10\\
\hline 
0 & 0 & 0 & 0 & 0 & 0 & 0 & 0 & 1 & 1 & 1\\
10 & 1 & 1 & 1 & 1 & 1 & 2 & 2 & 2 & 2 & 2\\
20 & 2 & 2 & 2 & 3 & 3 & 3 & 3 & 3 & 3 & 3\\
30 & 3 & 4 & 4 & 4 & 4 & 4 & 4 & 4 & 4 & 5\\
40 & 5 & 5 & 5 & 5 & 5 & 5 & 5 & 6 & 6 & 6\\
50 & 6 & 6 & 6 & 6 & 6 & 7 & 7 & 7 & 7 & 7\\
60 & 7 & 7 & 7 & 8 & 8 & 8 & 8 & 8 & 8 & 8\\
70 & 8 & 10 & 10 & 10 & 10 & 10 & 10 & 10 & 10 & 11\\
80 & 11 & 11 & 11 & 11 & 11 & 11 & 11 & 12 & 12 & 12\\
90 & 12 & 12 & 12 & 12 & 12 & 13 & 13 & 13 & 13 & 13\\
100 & 13 & 13 & 13 & 14 & 14 & 14 & 14 & 14 & 14 & 14\\
110 & 14 & 15 & 15 & 15 & 15 & 15 & 15 & 15 & 15 & 16\\
120 & 16 & 16 & 16 & 16 & 16 & 16 & 16 & 17 & 17 & 17\\
130 & 17 & 17 & 17 & 17 & 17 & 18 & 18 & 18 & 18 & 18\\
140 & 18 & 18 & 18 & 20 & 20 & 20 & 20 & 20 & 20 & 20\\
150 & 20 & 21 & 21 & 21 & 21 & 21 & 21 & 21 & 21 & 22\\
160 & 22 & 22 & 22 & 22 & 22 & 22 & 22 & 23 & 23 & 23\\
170 & 23 & 23 & 23 & 23 & 23 & 24 & 24 & 24 & 24 & 24\\
180 & 24 & 24 & 24 & 25 & 25 & 25 & 25 & 25 & 25 & 25\\
190 & 25 & 26 & 26 & 26 & 26 & 26 & 26 & 26 & 26 & 28\\
\hline \hline
\end{tabular}
\caption{The initial values of $v(n), \, n=1, 2, \dots, 200$.}
\label{table-initvalues-vn}
\end{center}
\end{tiny}
\end{table}

\begin{table}[h]
\begin{tiny}
\begin{center}
\begin{tabular}{|c|ccccc ccccc|}
\hline \hline
$n$&1& & & & 5& & & & & 10\\
\hline 
0 & 1 & 2 & 3 & 5 & 6 & 7 & 8 & 12 & 14 & 15\\
10 & 16 & 18 & 19 & 20 & 21 & 27 & 28 & 30 & 31 & 33\\
20 & 34 & 35 & 36 & 40 & 42 & 43 & 46 & 48 & 49 & 50\\
30 & 51 & 59 & 60 & 61 & 62 & 66 &  67 & 68 & 69 & 73\\
40 & 74 & 75 & 76 & 78 & 80 & 81 & 82 & 88 & 90 & 92\\
50 & 93 & 95 & 96 & 99 & 100 & 104 & 105 & 106 & 107 & 109\\
60 & 110 & 111 & 113 & 125 & 126 & 127 & 128 & 130 & 131 & 132\\
70 & 133 & 141 & 142 & 143 & 145 & 147 & 148 & 149 & 150 & 156\\
80 & 161 & 162 & 163 & 165 & 166 & 167 & 168 & 172 &  173 & 175\\
90 & 176 & 178 & 179 & 180 & 181 & 189 & 190 & 192 & 194 & 198\\
100 & 199 & 200 &  201 & 205 & 206 & 207 & 208 & 214 & 215 & 216\\
110 & 217 & 223 & 224 & 225 & 226 & 228 &  230 & 231 & 232 & 236\\
120 & 238 & 239 & 240 & 242 & 245 & 247 & 248 & 264 & 265 & 266\\
130 &  267 & 269 & 270 & 271 & 274 & 278 & 279 & 280 & 281 & 283\\
140 & 284 & 285 & 286 & 298 &  299 & 300 & 302 & 304 & 305 & 307\\
150 & 308 & 312 & 314 & 315 & 316 & 318 & 319 & 320 &  321 & 329\\
160 & 330 & 335 & 336 & 338 & 339 & 340 & 341 & 345 & 347 & 348\\
170 & 350 & 352 &  353 & 354 & 356 & 362 & 363 & 364 & 365 & 369\\
180 & 370 & 371 & 372 & 376 & 377 & 378 &  379 & 381 & 384 & 385\\
190 & 386 & 398 & 399 & 400 & 401 & 405 & 406 & 408 & 409 & 417\\
\hline \hline
\end{tabular}
\caption{The initial values of $w(n), \, n=1, 2, \dots, 200$.}
\label{table-initvalues-wn}
\end{center}
\end{tiny}
\end{table}

If we look at the sequences $\left\{a(n)\right\}_{n \in \bN}$ and 
$\left\{h(n)\right\}_{n \in \bN}$ from the inverse perspective, we can define two more sequences.
Let $S_a(n)$ denote the smallest number $k \in \bN$, for which exactly $n$ nonisomorphic
abelian groups of order $k$ exist (\cite{SloaneA046056}). The first 60 elements of the
sequence $\left\{S_a(n)\right\}_{n \in \bN}$ are given in Table \ref{table-initvalues-san}.
Here $0$ denotes the case, where $S_a(n)$ does not exist ($n$ is not a product of partition
numbers). These indices $n$ are exactly multiples of primes in the sequence $s$
(\cite{SloaneA046064}).

\begin{table}[h]
\begin{tiny}
\begin{center}
\begin{tabular}{|c|ccccc ccccc|}
\hline \hline
$n$&1& & & & 5& & & & & 10\\
\hline 
0 & 1 & 4 & 8 & 36 & 16 & 72 & 32 & 900 & 216 & 144 \\
10 & 64 & 1800 & 0 & 288 & 128 & 44100 & 0 & 5400 & 0 & 3600 \\
20 & 864 & 256 & 0 & 88200 & 1296 & 0 & 27000 & 7200 & 0 & 512 \\
30 & 0 & 5336100 & 1728 & 0 & 2592 & 264600 & 0 & 0 & 0 & 176400\\ 
40 & 0 & 1024 & 0 & 2304 & 3456 & 0 & 0 & 10672200 & 7776 & 32400\\
50 & 0 & 0 & 0 & 1323000 & 5184 & 2048 & 0 & 0 & 0 & 4608 \\
\hline \hline
\end{tabular}
\caption{The initial values of $S_a(n), \, n=1, 2, \dots, 60$, (\cite{SloaneA046056}, A046056).}
\label{table-initvalues-san}
\end{center}
\end{tiny}
\end{table}

Let $S_h(n)$ denote the smallest number $k \in \bN$, for which exactly $n$ nonisomorphic
hamiltonian groups of order $k$ exist. The first 30 elements of the sequence 
$\left\{S_h(n)\right\}_{n \in \bN}$ are given in Table \ref{table-initvalues-shn}, where 
again $0$ denotes the case, where $n$ is not a product of partition numbers and $S_h(n)$ does
not exist. 

\begin{table}[h]
\begin{tiny}
\begin{center}
\begin{tabular}{|c|ccccc ccccc|}
\hline \hline
$n$&1& & & & 5& & & & & 10\\
\hline 
0 & 8 & 72 & 216 & 1800 & 648 & 5400 & 1944 & 88200 & 27000 & 16200\\
10 & 5832 & 264600 & 0 & 48600 & 17496 & 10672200 & 0 & 1323000 & 0 & 793800\\
20 & 243000 & 52488 & 0 & 32016600 & 405000 & 0 & 9261000 & 2381400 & 0 & 157464 \\
\hline \hline
\end{tabular}
\caption{The initial values of $S_h(n), \, n=1, 2, \dots, 30$.}
\label{table-initvalues-shn}
\end{center}
\end{tiny}
\end{table}

Let us finish with two open problems. Think of computing the genus of
each of the groups $\Gamma \in {\cal A} \cup {\cal H}$, counted by $b(n)$. 
Since $\bZ_n \in {\cal A} \cup {\cal H}$, the minimal genus is $0$.
A natural question is therefore to determine
$$g(n) := \max \{\gamma(\Gamma)|\,  \Gamma \in {\cal A} \cup {\cal H},\, 
|\Gamma | = n\}.$$
The sequence $\left\{ g(n)\right \}_{n \in \bN} = \left( 0,0,0,0,0,0,0,1,\dots \right).$

Another interesting problem is a generalization of the considered problem, namely,
the problem of determining the number of groups, whose every subgroup is 2-subnormal 
(\cite{Mah86, Sta78}). A subgroup $H$ of group $G$ is said to be \emph{2-subnormal} in $G$
if there is a series
$$H = H_2 \triangleleft H_1 \triangleleft H_0 = G$$
of subgroups in $G$ (see \cite{LenSto87}). 
Such a subgroup is said also to be \emph{of defect} 2. Similarly, subgroups $H$ of defect 1
in $G$ are precisely normal subgroups of $G$.

\section{Acknowledgements}

Research was supported in part by a grant J1-6062 from Ministrstvo za
\v{s}olstvo, znanost in \v{s}port Republike Slovenije. Part of the
research was conducted while the first author was visiting Neil R.
Grabois Professor of Mathematics at Colgate University.

\end{document}